\documentclass[12pt,a4paper]{article}
\usepackage{amssymb}

\usepackage{amssymb}
\usepackage{hyperref}

\usepackage{graphicx,amssymb,amsfonts,epsfig,a4,amsmath,amsthm,xypic}
\usepackage[latin1]{inputenc}

\usepackage{tikz-cd}

\newtheorem{Thm}{Theorem}[section]
\newtheorem{Prop}[Thm]{Proposition}
\newtheorem{Lem}[Thm]{Lemma}
\newtheorem{Cor}[Thm]{Corollary}

\newtheorem{Rem}[Thm]{Remark}
\newtheorem{Ex}{Example}

\newcommand{\Z}{\mathbb{Z}}

\newcommand{\R}{\mathbb{R}}
\newcommand{\C}{\mathbb{C}}
\newcommand{\Q}{\mathbb{Q}}

\newcommand{\T}{\mathbb{T}}

\newcommand{\bpr}{\noindent \textbf{Proof}: }

\newcommand{\epr}{~$\blacksquare$}

\newcommand{\ra}{\rightarrow}
\newcommand{\tr}{\operatorname{tr}}
\title{Rieffel projections and 2-by-2 matrices}
\author{Olivier ISELY and Alain VALETTE\footnote{The second author acknowledges support of the Institut Henri Poincar\'e (UAR 839 CNRS-Sorbonne Universit\'e), and LabEx CARMIN (ANR-10-LABX-59-01).}}
\date{May 2, 2025}

\begin{document}

\maketitle

\begin{abstract} For a compact space $Y$, we view $C(Y\times S^1)$ as the crossed product $C(Y)\rtimes\Z$, with $\Z$ acting trivially. This allows us to study Rieffel projections in $M_2(C(Y\times S^1)$: we characterize them and compute their image under the projection $\partial_0:K_0(C(Y\times S^1))\rightarrow K_1(C(Y))$. We provide a new Rieffel projection in $M_2(C(\T^2))$, different from Loring's one \cite{Lor}, and involving only trigonometric polynomials plus the square root of $2-e^{2\pi i\theta}-e^{-2\pi i\theta}$. We give applications of this projection, e.g. explicit generators for the K-theory of $C(\T^3)$. Finally, we prove that, if a Banach algebra completion $\mathcal{B}$ of $\C[\Z^n]$ is continuously contained in $C(\T^n)$ and such that the Fourier series of $(2-e^{2\pi i\theta_j}-e^{-2\pi i\theta_j})^{1/2}\;(j=1,...,n)$ converges in $\mathcal{B}$, then the inclusion $\mathcal{B}\hookrightarrow C(\T^n)$ induces isomorphisms in K-theory.
\end{abstract}

\section{Introduction}

Let $\alpha$ be a *-automorphism of a unital $C^*$-algebra $A$. The crossed product $A\rtimes_\alpha\Z$ is the universal $C^*$-algebra generated by $A$ and a unitary $u$ such that $u^*au=\alpha(a)$ for every $a\in A$. A {\it Rieffel projection} is a projection in $A\rtimes_\alpha\Z$ of the form $p=u^*x_1^*+x_0+x_1u$, with $x_0,x_1\in A$, and $x_0^*=x_0$. The conditions $p=p^*=p^2$ are equivalent to
$$0=\alpha(x_1)x_1;$$
$$x_1=x_0x_1+x_1\alpha^{-1}(x_0);$$
$$x_0=x_0^2+x_1x_1^*+\alpha(x_1^*x_1).$$
Rieffel projections were introduced in \cite{Rie} to show that the irrational rotation algebra $A_\theta$, viewed as the crossed product of $C(S^1)$ by the rotation of angle $2\pi\theta\;(\theta\notin\Q)$, contains non-trivial projections.

The celebrated Pimsner-Voiculescu 6-terms exact sequence \cite{PV80} connects the K-theory of $A\rtimes_\alpha \Z$ with the K-theory of $A$:
$$\begin{array}{ccccc}
K_0(A) & \stackrel{1-\alpha_*}{\longrightarrow} & K_0(A) & \stackrel{\iota_*}{\longrightarrow}& K_0(A\rtimes_\alpha \Z) \\
& & & & \\
\partial_1 \uparrow &  &  &  & \downarrow \partial_0 \\
& & & & \\
K_1(A\rtimes_\alpha \Z) & \stackrel{\iota_*}{\longleftarrow} & K_1(A)) & \stackrel{1-\alpha_*}{\longleftarrow} & K_1(A)
\end{array}$$
where $\iota:A\rightarrow A\rtimes_\alpha \Z$ is the canonical inclusion.

In the Appendix of \cite{PV80}, Pimsner and Voiculescu explain how to compute the image of the K-theory class of a Rieffel projection under the map $\partial_0$:

\begin{Thm}\label{PV} Let $p=u^*x_1^*+x_0+x_1u$ be a Rieffel projection in $A\rtimes_\alpha \Z$. Denote by $\ell_{x_1}$ the left support projection of $x_1$ in the enveloping von Neumann algebra of $A$. Then $x_0\ell_{x_1}$ is self-adjoint (equivalently $[x_0,\ell_{x_1}]=0$), the unitary $\exp(2\pi ix_0\ell_{x_1})$ belongs to $A$, and
$$\partial_0[p]=[\exp(2\pi ix_0\ell_{x_1})].$$\epr
\end{Thm}

When $\Z$ is acting trivially on $A$, so that $A\rtimes_{Id}\Z=A\otimes C(S^1)$, the Pimsner-Voiculescu exact sequence degenerates into two short exact sequences:
$$0\rightarrow K_i(A)\stackrel{\iota_*}{\rightarrow}K_i(A\otimes C(S^1))\stackrel{\partial_i}{\rightarrow}K_{1-i}(A)\rightarrow 0\;(i=0,1).$$
Let $\varepsilon:A\otimes C(S^1)\rightarrow A$ denote evaluation at $1\in S^1$; since $\varepsilon$ provides a left inverse for $\iota$, we recover the classical isomorphisms
$$K_i(A\otimes C(S^1))=K_0(A)\oplus K_1(A)\;(i=0,1),$$
with $\partial_i$ becoming the projection $K_i(A\otimes C(S^1))\rightarrow K_{1-i}(A)$ in the above decomposition.

We focus on $i=0$. It turns out that Rieffel projections can be interesting even when $\Z$ is acting trivially! That was first observed by Loring \cite{Lor}: if $\T^2$ denotes the 2-torus, it is well-known that $K_0(C(\T^2))=\Z^2$, with one generator being the class $[1]$ of the unit in $C(\T^2)$, the other generator, sometimes called the Bott element, being non-trivial. Viewing $C(\T^2)$ as the crossed product $C(S^1)\rtimes_{Id}\Z$, Loring was able to express the non-trivial generator as the class of a Rieffel projection in $M_2(C(\T^2))$, therefore the projection on a rank 1 sub-bundle of the trivial bundle $\T^2\times\C^2$.

Motivated by Loring's paper \cite{Lor}, for $Y$ a compact space, we view $C(Y\times S^1)$ as the crossed product $C(Y)\rtimes_{Id}\Z$ and we consider the question of finding Rieffel projections in $M_2(C(Y\times S^1))$, i.e. projections of the form 
\begin{equation}\label{Rieffel}
P(y,\varphi)=\exp(-2\pi i\varphi)X_1(y)^* + X_0(y) + \exp(2\pi i\varphi)X_1(y),
\end{equation}
 with $\varphi\in[0,1], X_0,X_1\in M_2(C(Y))$, and $X_0=X_0^*$. We will say that such a Rieffel projection $P$ is {\it non-trivial} if $X_1\neq 0$, and that it has {\it rank 1} if $P(y,\varphi)$ has rank 1 for every $(y,\varphi)\in Y\times [0,1]$ (note that the rank 1 condition follows from non-triviality when $Y$ is connected). Non-trivial, rank 1 Rieffel projections are characterized as follows (see section 2 for the proof):
 
 \begin{Prop}\label{main1} There exists a non-trivial, rank 1, Rieffel projection in \linebreak $M_2(C(Y\times S^1))$ if and only if there exists $a,b,\alpha,\beta,\gamma\in C(Y)$, with $a$ real-valued, $|\beta|^2+|\gamma|^2\neq 0$, such that:
 \begin{itemize}
 \item $\beta\gamma=-\alpha^2;$
 \item $(2a-1)\alpha+b\gamma+\overline{b}\beta=0;$
 \item $a=a^2+|b|^2+2|\alpha|^2+|\beta|^2+|\gamma|^2.$
 \end{itemize}
 In that case $X_0=\left(\begin{array}{cc}a & b \\ \overline{b} & 1-a\end{array}\right)$ and $X_1=\left(\begin{array}{cc}\alpha & \beta \\\gamma & -\alpha\end{array}\right)$.
 \end{Prop}
 
 \begin{Ex}\label{trivial1} Let $a,\beta$ in $C(Y)$ be such that $a$ is real-valued. The self-adjoint matrix $\left(\begin{array}{cc}a & \beta \\\overline{\beta} & 1-a\end{array}\right)$ is a rank 1 projection in $M_2(C(Y))$ if and only if $a=a^2+|\beta|^2$. In that case, 
 $$P(y,\varphi)=\left(\begin{array}{cc}a(y) & \beta(y)\exp(2\pi i\varphi) \\\overline{\beta(y)}\exp(-2\pi i\varphi) & 1-a(y)\end{array}\right)$$
 is a projection in $M_2(C(Y\times S^1))$, which is indeed a Rieffel projection with $X_0=\left(\begin{array}{cc}a & 0 \\0 & 1-a\end{array}\right)$ and $X_1=\left(\begin{array}{cc}0 & \beta \\0 & 0\end{array}\right)$. We call this Rieffel projection of {\bf trivial type}: indeed we will see in Example \ref{Loring1} below, that its K-theory class is always in the kernel of $\partial_0:K_0(C(Y\times S^1))\rightarrow K_1(C(Y))$.
 \end{Ex}
 
To identify the image of the class of a general Rieffel projection in \linebreak $M_2(C(Y\times S^1))$ under $\partial_0:K_0(C(Y\times S^1))\rightarrow K_1(C(Y))$, we need two more notations. For $y\in Y$, let $\ell_{X_1}(y)$ denote the projection on the range of $X_1(y)$, so that $\ell_{X_1}$ is the left support projection of $X_1$ in the enveloping von Neumann algebra of $M_2(C(Y))$, as in Theorem \ref{PV}. Set $supp(X_1)=\{y\in Y:X_1(y)\neq 0\}$. With the help of Theorem \ref{PV}, the next result will be proved in Section 3.

\begin{Thm}\label{main2} Let $P(y,\varphi)=\exp(-2\pi i\varphi)X_1(y)^* + X_0(y) + \exp(2\pi i\varphi)X_1(y)$ be a non-trivial rank 1 Rieffel projection in $M_2(C(Y\times S^1))$. There exists a function $f:Y\rightarrow [0,1]$, continuous on $supp(X_1)$, vanishing elsewhere, such that:
\begin{itemize}
\item $X_0\ell_{X_1}=f\ell_{X_1};$
\item $\exp(2\pi if)\in C(Y);$
\item $\partial_0[P]$ is the class of the unitary $(1-\ell_{X_1})+\exp(2\pi if)\ell_{X_1}$ in $K_1(C(Y));$
\item if $supp(X_1)$ is contractible, then $\partial_0[P]=[\exp(2\pi if)]$ in $K_1(C(Y))$.
\end{itemize}
\end{Thm}
 
 \begin{Ex}\label{Loring1} We will say that a Rieffel projection $P(y,\varphi)\in M_2(C(Y\times S^1))$ is of {\bf Loring type}, if it is of the form 
 \begin{equation}\label{Loringguess}
 P(y,\varphi)=\left(\begin{array}{cc}a(y) & b(y)+\beta(y)\exp(2\pi i\varphi) \\b(y)+\beta(y)\exp(-2\pi i\varphi) & 1-a(y)\end{array}\right),
 \end{equation}
where $a,b,\beta$ are continuous, real-valued on $Y$, and such that $b\beta=0$ and $a=a^2 + b^2+\beta^2$ (cf Proposition \ref{main1}). 
Hence:
 $$X_0=\left(\begin{array}{cc}a & b \\b & 1-a\end{array}\right);\;\;X_1=\left(\begin{array}{cc}0 & \beta \\0 &0\end{array}\right).$$
 Then $\ell_{X_1}(y)=\left(\begin{array}{cc}1 & 0 \\0 & 0\end{array}\right)$ when $\beta(y)\neq 0$, while $\ell_{X_1}(y)=0$ when $\beta(y)=0$. And 
 $$f(y)=\left\{\begin{array}{ccc}a(y) & if & \beta(y)\neq 0 \\0 & if & \beta(y)=0\end{array}\right..$$
 So that, for a projection $P(y,\varphi)$ of Loring type, Theorem \ref{main2} says that $\partial_0[P]=[u]$, where $u$ is the unitary in $C(Y)$ given by:
 $$u(y)=\left\{\begin{array}{ccc}\exp(2\pi i a(y)) & if & \beta(y)\neq 0 \\1 & if & \beta(y)=0\end{array}\right..$$
 
 Projections of trivial type (see Example \ref{trivial1}) are projections $P$ of Loring type with $b=0$. From $a=a^2+\beta^2$ it follows that $a(y)$ is either 0 or 1 when $\beta(y)=0$. This means that $u(y)=\exp(2\pi ia(y))$ for every $y\in Y$. So $u$ admits a continuous branch of the logarithm on $Y$, hence $\partial_0[P]=[u]=0$ in $K_1(C(Y))$, as announced in Example \ref{trivial1}.
 \end{Ex}
 
\begin{Ex}\label{Loring2} In the case $Y=S^1$, the projections from Example \ref{Loring1} were introduced by Loring in section 2 of \cite{Lor}\footnote{Quoting from page 200 of \cite{Lor}: {\it ``Following along the lines of \cite{Rie}, we guess that the desired projection can be taken to be of the form (\ref{Loringguess})''}.}. Taking a continuous, 1-periodic function $a$ with values in $[0,1]$, decreasing from 1 to 0 on $[0,1/2]$ and increasing from 0 to 1 on $[1/2,1]$, and defining moreover:
$$b(\theta)=\left\{\begin{array}{ccc}\sqrt{a(\theta)-a(\theta)^2}  & if & \theta\in [0,1/2]\\ 0 &  if & \theta\in [1/2,1] \end{array}\right.$$
$$\beta(\theta)=\left\{\begin{array}{ccc}0 & if & \theta\in [0,1/2]\\ \sqrt{a(\theta)-a(\theta)^2}  & if & \theta\in [1/2,1] \end{array}\right.,$$
Loring gets a Rieffel projection $P_L(\theta,\varphi)$ defining the non-trivial generator of $K_0(C(\T^2))$, sometimes called the Bott element of $\T^2$; Loring's proof of this fact (Proposition 2.1 in \cite{Lor}) was in the spirit of non-commutative differential geometry; it also follows from Example \ref{Loring1}: indeed the unitary function $u$ is then given by
$$u(\theta)=\left\{\begin{array}{ccc}1 & if & \theta\in [0,1/2] \\\exp(2\pi ia(\theta)) & if & \theta\in [1/2,1]\end{array}\right..$$
Since $a$ increases from 0 to 1 on $[1/2,1]$, the function $u$ has winding number 1, hence $\partial_0[P_L]=[\exp(2\pi i\theta)]$, and $P_L$ defines the Bott element of $\T^2$; our proof is more in the spirit of non-commutative topology.
\end{Ex}

As noticed by Loring, via Fourier series we may identify $C(\T^2)$ with the group $C^*$-algebra $C^*(\Z^2)$, i.e. the universal $C^*$-algebra generated by two commuting unitaries, so any pair $U,V$ of commuting unitaries in a unital $C^*$-algebra $A$ defines a *-homomorphism $C(\T^2)\rightarrow C^*(U,V)$ which is onto. By functional calculus \`a la Gelfand the projection $P_L$ of Example \ref{Loring2} can be applied to $U,V$, and defines the projection 
$$P_L(U,V)=\left(\begin{array}{cc}a(U) & b(U)+\beta(U)V \\b(U)+\beta(U)V^* & 1-a(U)\end{array}\right)$$
in $M_2(C^*(\Z^2))$, with $a,b,\beta$ as in Example \ref{Loring2}.

In view of the definitions of $a,b,\beta$, it is clear that the full force of Gelfand functional calculus is used to define $a(U), b(U), \beta(U)$. It may be asked whether there are other formulae, more algebraic, also defining the Bott element of $\T^2$. We first observe that such a formula cannot be completely algebraic, i.e. it cannot involve trigonometric polynomials only. The reason is that, by a result by Bass, Heller and Swan (see the Corollary on the first page of \cite{BHS}), for the ring $\C[t_1,...,t_n,t_1^{-1},...,t_n^{-1}]=\C[\Z^n]$ of Laurent series with complex coefficients, the inclusion 
of $\C$ into $\C[\Z^n]$ induces an isomorphism in $K_0$, so the equality $K_0(\C[\Z^2])=\Z$ precludes the existence of a purely polynomial formula. However, if we allow algebraic functions (in our case, square roots), such a formula does exist.

\begin{Thm}\label{magict} The following formula defines a projection in $C(\T^2)$:

\begin{equation}\label{magic0}
\begin{scriptsize}
P_I(\theta,\varphi)=\frac{1}{4}\left(\begin{array}{cc}4-(1-\cos(2\pi\theta))(1-\cos(2\pi\varphi)) &  \sin(2\pi\theta)(1-\cos(2\pi\varphi))+i\sin(2\pi\varphi)\sqrt{2(1-\cos(2\pi\theta))}
\\ \sin(2\pi\theta)(1-\cos(2\pi\varphi))-i\sin(2\pi\varphi)\sqrt{2(1-\cos(2\pi\theta))} & (1-\cos(2\pi\theta))(1-\cos(2\pi\varphi)) \end{array}\right).
\end{scriptsize}
\end{equation}
It is actually a Rieffel projection associated with
\begin{equation}\label{magic1}
X_0(\theta)=\frac{1}{4}\left(\begin{array}{cc}3+\cos(2\pi\theta) & \sin(2\pi\theta)\\\sin(2\pi\theta) & 1-\cos(2\pi\theta) \end{array}\right);
\end{equation}
\begin{equation}\label{magic2}
X_1(\theta)=\frac{1}{8}\left(\begin{array}{cc}1-\cos(2\pi\theta) & \sqrt{2(1-\cos(2\pi\theta))}-\sin(2\pi\theta) \\-\sqrt{2(1-\cos(2\pi\theta))}-\sin(2\pi\theta) & -1+\cos(2\pi\theta)\end{array}\right).
\end{equation}
Moreover $\partial_0[P_I]=[\exp(2\pi i\theta)]$, hence $P_I$ defines the Bott element of $\T^2$.
\end{Thm}

Formula (\ref{magic0}) follows from formulae (\ref{magic1}), (\ref{magic2}) by direct computation; the latter formulae appeared in the first author's PhD thesis (see section 3.1 in \cite{Isely}). In section 4 we provide an explanation for them.

\medskip
Applying $P_I$ to a pair $U,V$ of commuting unitaries in a $C^*$-algebra $A$, yields the projection

\begin{scriptsize}
\begin{equation}\label{magicproj}
P_I(U,V)=\frac{1}{16}\left(\begin{array}{cc}16-(2-U-U^*)(2-V-V^*) & 2(2-U-U^*)^{1/2}(V-V^*)+i(U^*-U)(2-V-V^*) \\ 2(2-U-U^*)^{1/2}(V^*-V)+i(U-U^*)(2-V-V^*)& (2-U-U^*)(2-V-V^*)\end{array}\right)
\end{equation}
\end{scriptsize}
in $M_2(A)$. Non-triviality of this projection in K-theory, is ensured by our Corollary \ref{crossedproduct}: if $\alpha$ is a *-automorphism of the unital $C^*$-algebra $A$, assume that $\alpha$ fixes some unitary $v\in A$. Then, with $u$ the canonical unitary implementing $\alpha$ in the crossed product $A\rtimes_\alpha\Z$, we have $\partial_0[P_I(v,u)]=[v]$ in $K_1(A)$. Other applications of Theorem \ref{magict} discussed in Section 5 are explicit degree 1 maps from $\T^2$ to the complex projective line $P^1(\C)$ (observe indeed that, from formula (\ref{magic0}), we have $P_I(\theta,0)=P_I(0,\phi)=\left(\begin{array}{cc}1 & 0 \\0 & 0\end{array}\right)$ for every $\theta,\phi\in [0,1]$, so that $P_I$ factors through the 2-sphere $S^2$); and explicit generators for $K_0$ and $K_1$ of the 3-torus $\T^3$. 

The final section 6 has a different flavour: as we noticed already, in (\ref{magic0}) we use only trigonometric polynomials, plus the square root of the simple polynomial $2(1-\cos(2\pi\theta))=2-e^{2\pi i\theta}-e^{-2\pi i\theta}$. Recall the Fourier series of $(2-e^{2\pi i\theta}-e^{-2\pi i\theta})^{1/2}$ (see e.g. formula 24.11 in \cite{SLL}):
\begin{equation}\label{fourier}
(2-e^{2\pi i\theta}-e^{-2\pi i\theta})^{1/2}=|\sin(\pi\theta)|=\frac{2}{\pi}-\frac{4}{\pi}\sum_{k=1}^\infty\frac{\cos(2k\pi\theta)}{4k^2-1}.
\end{equation}

We prove that, if a Banach algebra completion of the group ring $\C[\Z^n]$ on the one hand embeds continuously in $C^*(\Z^n)=C(\T^n)$, on the other hand, for $s_i$ in the canonical basis of $\Z^n$, contains the square roots of the group ring elements $2-s_i-s_i^*$ (in the sense that the series (\ref{fourier}) converges in the algebra),  then the algebra has the same K-theory as $C(\T^n)$. Precisely:

\begin{Prop}\label{sameKth} Fix $n\geq 1$. Let $\mathcal{B}(\Z^n)$ be a Banach algebra completion of $\C[\Z^n]$, continuously contained in  $C^*(\Z^n)=C(\T^n)$. Let $\{s_1,...,s_n\}$ be the canonical basis of $\Z^n$. Assume that, for every $i=1,...,n$, the series 
\begin{equation}\label{root}
\frac{2}{\pi}(1-\sum_{k=1}^\infty \frac{s_i^k+(s_i^*)^k}{4k^2-1}),
\end{equation}
(defining $(2-s_i-s_i^*)^{1/2}$) converges in $\mathcal{B}(\Z^n)$. Then the inclusion $\mathcal{B}(\Z^n)\hookrightarrow C^*(\Z^n)$ induces isomorphism in K-theory.
\end{Prop}

In lemma \ref{spectral} we actually prove the stronger statement that $\mathcal{B}(\Z^n)$ is spectral in $C^*(\Z^n)$, i.e. an element of $\mathcal{B}(\Z^n)$ is invertible in $C^*(\Z^n)$ if and only if it is invertible in $\mathcal{B}(\Z^n)$. 
Observe that the algebra $\mathcal{B}(\Z^n)=\ell^1(\Z^n)$ satisfies the assumptions of Proposition \ref{sameKth}. 

In the concluding remarks, we summarize the known results on embeddings of Banach algebras inducing isomorphisms in K-theory, in particular by J.-B. Bost \cite{JBB} on the Oka principle, and by V. Lafforgue \cite{Laff} on the Baum-Connes conjecture.

{\bf Acknowledgements:} We thank M.-P. Gomez-Aparicio for pointing out references \cite{Nica,SaWi}, N. Higson for mentioning reference \cite{BHS}, P. Julg for reminding us of Wiener's $1/f$-theorem (this remark eventually led to Proposition \ref{sameKth}), and G. Skandalis for detecting a mistake in the original proof of lemma \ref{spectral}.

\section{Proof of Proposition \ref{main1}}

We start with a lemma from linear algebra. We denote by ${\bf 1}_n$ the $n$-by-$n$ identity matrix.

\begin{Lem}\label{linalg} Let $X\in M_2(\C)$. \begin{enumerate}
\item The matrix $X$ is nilpotent (i.e. $X^n=0$ for some $n\geq 2$) if and only if $X^2=0$, if and only if $\tr(X)=\det(X)=0$.
\item The matrix $X$ is a rank 1 idempotent if and only if $\tr(X)=1$ and $\det(X)=\nolinebreak 0$.
\item If $X$ is nilpotent, then $X^*X+XX^*=\tr(X^*X){\bf 1}_2.$
\end{enumerate}\end{Lem}

\bpr \begin{enumerate}
\item By Cayley-Hamilton, the characteristic polynomial of $X$, namely \linebreak $t^2-\tr(X)t+\det(X)$, vanishes on $X$. To prove the statement, clearly we may assume $X\neq 0$: so the polynomial $t^n$ also vanishes on $X$, which happens if and only if the polynomial $t^2$ vanishes on $X$, if and only if $\tr(X)=\det(X)=0$. 
\item The proof is very similar to the one of 1.
\item This can be proved by direct computation, let us give a more conceptual proof. We may assume $X\neq 0$. Then $X^*X$ is a positive multiple of the projection $P$ on $Im(X)=\ker(X)$, while $XX^*$ is a positive multiple of the projection $Q$ on $Im(X^*)=(\ker(X))^\perp$. Taking traces, we see that those multiples coincide and are equal to $\tr(X^*X)=\tr(XX^*)$. Finally $P+Q={\bf 1}_2$.\epr
\end{enumerate}

To prove Proposition \ref{main1}, we look for Rieffel projections in $M_2(C(Y\times S^1))$ of the form (\ref{Rieffel}), which corresponds to a pair of matrices $X_0,X_1\in M_2(C(Y))$ with $X_0=X_0^*$ and
\begin{equation}\label{Rieffel1}
X_1^2=0;
\end{equation}
\begin{equation}\label{Rieffel2}
X_1=X_0X_1+X_1X_0;
\end{equation}
\begin{equation}\label{Rieffel3}
X_0=X_0^2 + X_1^*X_1+X_1X_1^*=X_0^2+\tr(X_1^*X_1){\bf 1}_2
\end{equation}
(where the last equality follows from lemma \ref{linalg}).
Moreover we want this Rieffel projection to be non-trivial (i.e. $X_1\neq 0$) and of rank 1, so that $\tr(X_0(y))=1$ for every $y\in Y$. Assuming that such a pair $X_0,X_1$ exists, we may write 
$$X_0=\left(\begin{array}{cc}a & b \\ \overline{b} & 1-a\end{array}\right);\;\;\;X_1=\left(\begin{array}{cc}\alpha & \beta \\\gamma & -\alpha\end{array}\right)$$
with $a,b,\alpha,\beta,\gamma\in C(Y)$ and $a$ real-valued. By direct computation, one checks that the condition (\ref{Rieffel1}) translates into 
\begin{equation}\label{nonlincond}
\beta\gamma=-\alpha^2
\end{equation}(see lemma \ref{linalg}); the condition (\ref{Rieffel2}) translates into 
\begin{equation}\label{linearcond}
(2a-1)\alpha+b\gamma+\overline{b}\beta=0;
\end{equation} and the condition (\ref{Rieffel3}) translates into 
\begin{equation}
a=a^2+|b|^2+2|\alpha|^2+|\beta|^2+|\gamma|^2.
\end{equation} This concludes the proof of Proposition \ref{main1}.
\epr

\section{Proof of Theorem \ref{main2}}

We denote by $Proj_1(M_n(\C))$ the space of rank 1 projections in $M_n(\C)$ (it is canonically homeomorphic to the $(n-1)$-dimensional complex projective space $P^{n-1}(\C)$).

\begin{Lem}\label{homotopy} Let $U$ be an open subset in the compact space $Y$, and let $f:Y\rightarrow\R$ be a function which is continuous on $U$, and 0 elsewhere.
\begin{enumerate}
\item The following are equivalent:
\begin{enumerate}
\item The function $\exp(2\pi if)$ is continuous on $Y$;
\item For any function $P:Y\rightarrow Proj_1(M_n(\C))$ which is continuous on $U$ and zero elsewhere, the map $\exp(2\pi i fP)$ is continuous from $Y$ to $U_n(\C)$.
\end{enumerate}
\item If $f,P$ are as above (with $\exp(2\pi if)$ continuous on $Y$), the K-theory class of $\exp(2\pi i fP)$ in $K_1(C(Y))$ only depends on the homotopy class of $P|_U$ as a map $U\rightarrow Proj_1(M_n(\C))$.
\item If moreover $U$ is contractible, then $[\exp(2\pi i fP)]=[\exp(2\pi i f)]$ in $K_1(C(Y))$.
\end{enumerate}
\end{Lem}

\bpr \begin{enumerate}
\item For $(b)\Rightarrow(a)$, assume that $\exp(2\pi i fP)$ is a continuous map on $Y$. So, applying the determinant, $\det(\exp(2\pi i fP)=\exp(\tr(2\pi ifP))$ is a continuous function on $Y$. Using $\tr(P)=1$ on $U$, we see that the equality $\exp(\tr(2\pi ifP))=\exp(2\pi if)$ holds everywhere on $Y$, hence the continuity of $\exp(2\pi if)$.

For $(a)\Rightarrow (b)$, suppose that $\exp(2\pi if)$ is continuous on $Y$. Using the assumption on $P$, the map $\exp(2\pi ifP)$ is certainly continuous on $U$. Let us show that it is also continuous at any $y\in Y\setminus U$. So let $(y_j)_{j\in I}$ be a net in $Y$ converging to $y$. By the Taylor expansion of the exponential:
\begin{equation}\label{Taylor}
\exp(2\pi ifP)=\exp(2\pi if)P + ({\bf 1}_n-P).
\end{equation}
Hence, denoting by $\|\cdot\|$ the operator norm on $M_n(\C)$, we have (using (\ref{Taylor}) for the second equality):
$$\|\exp(2\pi i f(y)P(y))-\exp(2\pi if(y_j)P(y_j))\|=\|{\bf 1}_n-\exp(2\pi if(y_j)P(y_j))\|$$
$$=\|(1-\exp(2\pi if(y_j)))P(y_j)\|\leq |\exp(2\pi if(y))-\exp(2\pi if(y_j))|$$
which goes to 0 for $j\rightarrow\infty$.

\item Let $P_0,P_1:Y\rightarrow Proj_1(M_n(\C))$ be maps, vanishing on $Y\setminus U$, such that $P_0|_U,P_1|_U$ are homotopic. Let $P:U\times [0,1]\rightarrow Proj_1(M_n(\C))$ be a homotopy, i.e. $P(y,0)=P_0(y)$ and $P(y,1)=P_1(y)$ for every $y\in U$. Extend $P$ to $Y\times[0,1]$ by 0 outside of $U\times [0,1]$. View $f$ as a function on $Y\times [0,1]$ by composing with the first projection. Then apply the first part of the lemma to $f$ and $P$, with $Y$ replaced by $Y\times [0,1]$ and $U$ replaced by $U\times [0,1]$. So $\exp(2\pi ifP)$ is continuous on $Y\times[0,1]$, and yields a homotopy between $\exp(2\pi ifP_0)$ and $\exp(2\pi ifP_1)$.

\item If $U$ is contractible, then $P$ is homotopic to $P_0$ which, on $U$, is a constant function to $Proj_1(M_n(\C))$, and is 0 on $Y\setminus U$. By the second part of the lemma together with (\ref{Taylor}), we have:
$$[\exp(2\pi ifP)]=[\exp(2\pi ifP_0)]=[\exp(2\pi if)P_0+(1-P_0)]=[\exp(2\pi if)].$$\epr
\end{enumerate}

{\bf Proof of Theorem \ref{main2}:}
Assume $X_1(y)\neq 0$. Then $\ell_{X_1}(y)$ is a rank 1 projection, and since $[X_0,\ell_{X_1}]=0$ (see Theorem \ref{PV}), $(\ell_{X_1}X_0)(y)$ is a self-adjoint matrix with image contained in the image of $\ell_{X_1}(y)$. So there exists a continuous function $f:supp(X_1)\rightarrow \R$ such that $(\ell_{X_1}X_0)(y)=f(y)\ell_{X_1}(y)$ for every \linebreak $y\in supp(X_1)$. We extend $f$ to $Y$ by setting it equal to 0 outside of $supp(X_1)$. Then, using (\ref{Taylor}):
$$\exp(2\pi iX_0\ell_{X_1})=\exp(2\pi if\ell_{X_1})=\exp(2\pi if)\ell_{X_1}+(1-\ell_{X_1}).$$
By Theorem \ref{PV}, $\exp(2\pi if)\ell_{X_1}+(1-\ell_{X_1})$ belongs to $M_2(C(Y))$; by lemma \ref{homotopy}, this is equivalent to $\exp(2\pi if)\in C(Y)$; finally, again by Theorem \ref{PV}, the image of the corresponding Rieffel projection under $\partial_0$, is the class of that unitary in $K_1(C(Y))$. If $supp(X_1)$ is contractible, by the 3rd part of lemma \ref{homotopy}, this coincides with the class of $\exp(2\pi if)$ in $C(Y)$.
\epr

It might be interesting to have explicit formulae, so write $X_0=\left(\begin{array}{cc}a & b \\ \overline{b} & 1-a\end{array}\right)$ and $X_1=\left(\begin{array}{cc}\alpha & \beta \\\gamma & -\alpha\end{array}\right)$ for some $a, b, \alpha,\beta,\gamma\in C(Y)$ (with $a$ real-valued), and observe that $supp(X_1)=supp(\beta)\cup supp(\gamma)$. 
\begin{itemize}
\item If $\beta(y)\neq 0$, then the vector $(\beta(y),-\alpha(y))$ spans the image of $X_1(y)$, hence 
$$\ell_{X_1}(y)=\frac{1}{|\alpha(y)|^2+|\beta(y)|^2}\left(\begin{array}{cc}|\beta(y)|^2 & -\beta(y)\overline{\alpha(y)} \\ -\overline{\beta(y)}\alpha(y)  & |\alpha(y)|^2\end{array}\right).$$
The (1,1)-coefficient of $(X_0\ell_{X_1})(y)$ is $\frac{-\beta(y)\overline{b(y)\alpha(y)}+ a(y)|\beta(y)|^2}{|\alpha(y)|^2+|\beta(y)|^2}$, hence 
$$f(y)= a(y)-\frac{\alpha(y)b(y)}{\beta(y)}.$$
For a projection of Loring type (see Example \ref{Loring1}), we have $\alpha=\gamma=0$, so for $\beta(y)\neq 0$ we find $f(y)=a(y)$, as in Example \ref{Loring1}.

\item If $\gamma(y)\neq 0$, then the vector $(\alpha(y),\gamma(y))$ spans the image of $X_1(y)$, hence 
$$\ell_{X_1}(y)=\frac{1}{|\alpha(y)|^2+|\gamma(y)|^2}\left(\begin{array}{cc}|\alpha(y)|^2 & \alpha(y)\overline{\gamma(y)} \\\overline{\alpha(y)}\gamma(y)  & |\gamma(y)|^2\end{array}\right).$$
The (2,2)-coefficient of $(X_0\ell_{X_1})(y)$ is $\frac{\alpha(y)\overline{b(y)\gamma(y)}+(1-a(y))|\gamma(y)|^2}{|\alpha(y)|^2+|\gamma(y)|^2}$, hence 
$$f(y)=1-a(y)+\frac{\alpha(y)\overline{b(y)}}{\gamma(y)}.$$
If $y\in supp(\beta)\cap supp(\gamma)$, then the two formulae for $f(y)$ match, as an easy consequence of (\ref{nonlincond}) and (\ref{linearcond}).
\end{itemize}

\section{Another Rieffel projection in $M_2(C(\T^2))$}

Here we exploit the action of the unitary group of $M_2(C(Y))$ on solutions $(X_0,X_1)$ of the system (\ref{Rieffel1}), (\ref{Rieffel2}), (\ref{Rieffel3}): for $U\in U_2(C(Y))$, the pair $(UX_0U^*,UX_1U^*)$ is a solution of the same system. Moreover, the function $f$ associated with the pair is invariant under this action, as follows from $X_0\ell_{X_1}=f\ell_{X_1}$ (see Theorem \ref{main2}).

For $Y=[0,1]$, let us start with the pair 
$$\mathcal{X}_0(\theta)=\left(\begin{array}{cc}\sin^2(\frac{\pi\theta}{2}) & 0 \\0 & \cos^2(\frac{\pi\theta}{2})  \end{array}\right),\;\mathcal{X}_1(\theta)=\left(\begin{array}{cc}0 & \sin(\frac{\pi\theta}{2})  \cos(\frac{\pi\theta}{2})  \\0 & 0\end{array}\right);$$
equations (\ref{Rieffel1}), (\ref{Rieffel2}), (\ref{Rieffel3}) are clearly satisfied, so we have a Rieffel projection - actually of trivial type - in $M_2(C(Y))$; here $\ell_{\mathcal{X}_1}(\theta)=\left(\begin{array}{cc}1 & 0 \\0 & 0\end{array}\right)$ for $0<\theta<1$, while $\ell_{\mathcal{X}_1}(0)=\ell_{\mathcal{X}_1}(1)=0$. Moreover 
\begin{equation}\label{goodf}
f(\theta)=\left\{\begin{array}{ccc}\sin^2(\frac{\pi\theta}{2}) & if & 0\leq\theta<1 \\0 & if & \theta=1 \end{array}\right.
\end{equation}

{\bf Proof of Theorem \ref{magict}:}
Observe that $\mathcal{X}_0(0)=\left(\begin{array}{cc}0 & 0 \\0 & 1\end{array}\right)$ and $\mathcal{X}_0(1)=\left(\begin{array}{cc}1 & 0 \\0 & 0\end{array}\right)$, so that $\mathcal{X}_0$ does not descend to $S^1$ (while $\mathcal{X}_1$ does). To correct this, we consider $U\in U_2(C([0,1]))$ given by a family of rotations:
$$U_\theta=\left(\begin{array}{cc}\sin(\frac{\pi\theta}{2}) & \cos(\frac{\pi\theta}{2}) \\-\cos(\frac{\pi\theta}{2}) & \sin(\frac{\pi\theta}{2})\end{array}\right);$$
$U_\theta$ is the rotation of angle $\frac{\pi}{2}(\theta-1),\;\theta\in [0,1]$.
Define then 
$$X_0(\theta)=U_\theta\mathcal{X}_0(\theta)U_\theta^*;\;X_1(\theta)=U_\theta\mathcal{X}_1(\theta)U_\theta^*,$$
which are now 1-periodic, as $X_0(0)=X_0(1)=\left(\begin{array}{cc}1 & 0 \\0 & 0\end{array}\right).$
Direct computations give:
$$X_0(\theta)=\left(\begin{array}{cc}\sin^4(\frac{\pi\theta}{2}) + \cos^4(\frac{\pi\theta}{2})  & \sin(\frac{\pi\theta}{2}) \cos^3(\frac{\pi\theta}{2}) -\sin^3(\frac{\pi\theta}{2}) \cos(\frac{\pi\theta}{2})  \\
\sin(\frac{\pi\theta}{2}) \cos^3(\frac{\pi\theta}{2}) -\sin^3(\frac{\pi\theta}{2}) \cos(\frac{\pi\theta}{2}) & 2\sin^2(\frac{\pi\theta}{2}) \cos^2(\frac{\pi\theta}{2}) \end{array}\right)$$
$$=\frac{1}{2}\left(\begin{array}{cc}2-\sin^2(\pi\theta) & \sin(\pi\theta) \cos(\pi\theta)  \\\sin(\pi\theta) \cos(\pi\theta)  & \sin^2(\pi\theta) \end{array}\right)$$
$$=\frac{1}{4}\left(\begin{array}{cc}3+\cos(2\pi\theta) & \sin(2\pi\theta)\\\sin(2\pi\theta) & 1-\cos(2\pi\theta) \end{array}\right).$$
Note that $X_0$ is indeed expressed in terms of 1-periodic trigonometric polynomials, i.e. it is in the algebra of 2-by-2 matrices over the complex group ring $\C[\Z]$.
$$X_1(\theta)=\sin(\frac{\pi\theta}{2})\cos(\frac{\pi\theta}{2})\left(\begin{array}{cc}\sin(\frac{\pi\theta}{2})\cos(\frac{\pi\theta}{2}) & \sin^2(\frac{\pi\theta}{2}) \\-\cos^2(\frac{\pi\theta}{2}) & -\sin(\frac{\pi\theta}{2})\cos(\frac{\pi\theta}{2})\end{array}\right)$$
$$=\frac{\sin(\pi\theta)}{4}\left(\begin{array}{cc}\sin(\pi\theta) & 1-\cos(\pi\theta) \\-1-\cos(\pi\theta) & -\sin(\pi\theta)\end{array}\right)$$
$$=\frac{1}{8}\left(\begin{array}{cc}1-\cos(2\pi\theta) & \sqrt{2(1-\cos(2\pi\theta))}-\sin(2\pi\theta) \\-\sqrt{2(1-\cos(2\pi\theta))}-\sin(2\pi\theta) & -1+\cos(2\pi\theta)\end{array}\right).$$
Note that $X_1$ is not in $M_2(\C[\Z])$ (it can't be, as mentioned in the Introduction).

It remains to prove that $\partial_0[P_I]=[\exp(2\pi i\theta)]$. But here $supp(X_1)=\,]0,1[$ is contractible, so lemma \ref{homotopy} applies, and $\partial_0[P_I]=[\exp(2\pi if)]$, where $f$ is given by (\ref{goodf}). But $\exp(2\pi if(\theta))=\exp(2\pi i\sin^2(\frac{\pi\theta}{2}))$ for every $\theta\in[0,1]$, and since $\exp(2\pi i\sin^2(\frac{\pi\theta}{2}))$ has winding number 1, the result follows. \epr

\section{Applications}

\subsection{Arbitrary crossed products with $\Z$}

Let $A$ be a unital $C^*$-algebra, and let $\alpha$ be any *-automorphism of $A$.

\begin{Cor}\label{crossedproduct} Assume that $\alpha(v)=v$ for some unitary $v\in A$. Then there exists a projection $p\in M_2(A\rtimes_\alpha \Z)$ such that $\partial_0[p]=[v]$ in $K_1(A)$.
\end{Cor}

\bpr Let $u$ be the unitary in $A\rtimes_\alpha \Z$ implementing the action. By assumption, $u$ and $v$ are commuting unitaries, so they define a surjective *-homomorphism $C(\T^2)\rightarrow C^*(u,v)$. Let then $p\in M_2(C^*(u,v))$ be defined either by $P_L(v,u)$ or by $P_I(v,u)$. By Example \ref{Loring2} or by Theorem \ref{magict}, we have $\partial_0[p]=[v]$ in $K_1(C^*(v))$. By naturality of the Pimsner-Voiculescu exact sequence (see Theorem 10.2.1 in \cite{Black}), we also have $\partial_0[p]=[v]$ in $K_1(A)$. \epr

\subsection{Maps of degree 1 from $\T^2$ to $P^1(\C)$}

Recall that $Proj_1(M_2(\C))$ denotes the space of rank 1 projections in $M_2(\C)$: it is canonically homeomorphic to the complex projective line $P^1(\C)$ (itself non-canonically homeomorphic to the 2-sphere $S^2$). A rank 1 projection in $M_2(C(Y))$ can be viewed as a continuous map $Y\rightarrow P^1(\C)$. From Example \ref{Loring2} and Theorem \ref{magict}, we immediately deduce:

\begin{Cor}\label{degree1} The maps $P_L$ and $P_I$ are degree 1 maps from $\T^2$ to $P^1(\C)$.
\epr
\end{Cor}

At this juncture, we observe that the Bass-Heller-Swan \cite{BHS} result already mentioned in the Introduction, implies that any map $\T^2\rightarrow Proj_1(M_2(\C))$ given by trigonometric polynomials, has degree 0. In the same vein, let us quote a result by Loday \cite{Loday}: viewing $\T^2$ as the product of two unit circles in $\R^2\times\R^2$, and $S^2$ as the unit sphere in $\R^3$, any polynomial map $\T^2\rightarrow S^2$ is homotopically trivial.

For more on degree 1 maps from $\T^2$ to $S^2$, see a discussion on Mathematics Stack Exchange \cite{MathStack}.

\subsection{K-theory of the 3-torus}

Let $\T^3$ denote the 3-torus. It is well-known that $K_0(C(\T^3))\simeq\Z^4\simeq K_1(C(\T^3))$. Writing $\T^3$ as the product of $\T^2$ and $S^1$ in various ways, we will exploit Theorem \ref{magict} to give explicit generators of $K_0(C(\T^3))$ and $K_1(C(\T^3))$. A point on $\T^3$ will be written $(\exp(2\pi i\theta),\exp(2\pi i\varphi),\exp(2\pi i\psi))$, with $\theta,\varphi,\psi\in [0,1]$.

\begin{Cor}\label{3torus}\begin{itemize}
\item The group $K_0(C(\T^3))$ is generated by the classes of the projections $1, P_I(\theta,\varphi),P_I(\varphi,\psi),P_I(\psi,\theta)$.
\item The group $K_1(C(\T^3))$ is generated by the classes of $\exp(2\pi i\theta),$ $\exp(2\pi i\varphi),$ $\exp(2\pi i\psi)$ and 
\begin{equation}\label{magicunit}
U_I(\theta,\varphi,\psi):=P_I(\theta,\varphi)\exp(2\pi i\psi)+({\bf 1}_2-P_I(\theta,\varphi))
\end{equation} (the latter being a unitary in $M_2(C(\T^3))$).
\end{itemize}\end{Cor}

\bpr\begin{itemize}
\item For $K_0$: for a compact space $Y$, we denote by $\widetilde{K}_0(C(Y))$ the {\it reduced} K-theory of $C(Y)$, namely $K_0(C(Y))/\Z\cdot[1]$. We define
$$\varepsilon_{\theta,\varphi}:C(\T^3)\rightarrow C(\T^2): f\mapsto [(\theta,\varphi)\mapsto f(\theta,\varphi,0)]$$
(evaluation at $\psi=0$); the maps $\varepsilon_{\varphi,\psi},\varepsilon_{\psi,\theta}$ are defined analogously. We consider then
$$\beta:=(\varepsilon_{\theta,\varphi})_*\oplus(\varepsilon_{\varphi,\psi})_*\oplus(\varepsilon_{\psi,\theta})_*:\widetilde{K}_0(C(\T^3))\rightarrow (\widetilde{K}_0(C(\T^2)))^3.$$
Note that both the source and the range of $\beta$ are isomorphic to $\Z^3$. We claim that $\beta[P_I(\theta,\varphi)]=([P_I(\theta,\varphi)],0,0)$. Indeed, clearly $\varepsilon_{\theta,\varphi}(P_I(\theta,\varphi))=P_I(\theta,\varphi)$, while $\varepsilon_{\varphi,\psi}(P_I(\theta,\varphi))=\left(\begin{array}{cc}1 & 0 \\0 & 0\end{array}\right)=\varepsilon_{\psi,\theta}(P_I(\theta,\varphi))$ by formula (\ref{magic0}); so the two latter vanish in reduced K-theory. Similarly we have $\beta[P_I(\varphi,\psi)]=(0,[P_I(\varphi,\psi)],0)$ and $\beta[P_I(\psi,\theta)]=(0,0,[P_I(\psi,\theta)])$, so that $\beta$ maps the set $\{[P_I(\theta,\varphi)],[P_I(\varphi,\psi)],[P_I(\psi,\theta)]\}\subset \widetilde{K}_0(C(\T^3))$ to the set
$$\{([P_I(\theta,\varphi)],0,0), (0,[P_I(\varphi,\psi)],0), (0,0,[P_I(\psi,\theta)])\}$$
which is a basis of $(\widetilde{K}_0(C(\T^2)))^3$ by Theorem \ref{magict}: this proves the desired result.
\item For $K_1$, we recall that, for $A,B$ unital $C^*$-algebras, there is the product map:
$$K_1(A)\otimes_\Z K_0(B)\rightarrow K_1(A\otimes B):[u]\otimes[p]\mapsto [u\otimes p+1\otimes(1-p)],$$
where $u\in U_m(A),p\in Proj(M_n(B))$ and $u\otimes p+1\otimes(1-p)\in U_{mn}(A\otimes B)$ (see e.g. Proposition 4.8.3 in \cite{HiRo}). By the K\"unneth formula (Theorem 23.1.3 in \cite{Black}), since all the involved K-groups are torsion-free, the product induces an isomorphism
$$K_1(C(T^2))\otimes_\Z K_0(C(S^1))\oplus K_0(C(\T^2))\otimes_\Z K_1(C(S^1))\simeq K_1(C(\T^3)).$$
The first summand contributes the generators $[\exp(2\pi i\theta)],[\exp(2\pi i\varphi)]$; \linebreak the second summand contributes the generators $[\exp(2\pi i\psi)]$ and \linebreak $[P_I(\theta,\varphi)\exp(2\pi i\psi)+({\bf 1}_2-P_I(\theta,\varphi))]$

\end{itemize}\epr

\begin{Rem} Corollary \ref{3torus} also holds with $P_I$ replaced all over by $P_L$. It may be interesting to notice that all generators of $K_i(C(\T^3))$ can be realized within $M_2(C(\T^3))$.
\end{Rem} 


\section{Spectral subalgebras of $C^*(\Z^n)$}

Let $A$ be a unital $C^*$-algebra, and $\mathcal{B}$ a unital subalgebra of $A$. Recall that $\mathcal{B}$ is said to be {\it spectral in $A$} if, for every element $x\in\mathcal{B}$ which is invertible in $A$, the inverse $x^{-1}$ belongs to $\mathcal{B}$. In other words, for any $x\in \mathcal{B}$, the spectrum $Sp_\mathcal{B}(x)$ of $x$ in $\mathcal{B}$ coincides with the spectrum $Sp_A(x)$ of $x$ in $A$ - hence the terminology.

Assume from now on that $\mathcal{B}$ is a Banach algebra continuously contained and dense in $A$. If $\mathcal{B}$ is spectral in $A$, then the inclusion $\mathcal{B}\hookrightarrow A$ induces isomorphisms in K-theory, e.g. by Theorem 1.3.1 in \cite{JBB}. So Proposition \ref{sameKth} will follow from:

\begin{Lem}\label{spectral} Under the assumptions of Proposition \ref{sameKth}, the subalgebra $\mathcal{B}$ is spectral in $C^*(\Z^n)=C(\T^n)$.
\end{Lem}

\bpr The proof is in two steps:
\begin{itemize}
\item We claim that $Sp_\mathcal{B}(s_i)=S^1$. Clearly $S^1=Sp_{C^*(\Z^n)}(s_i)$ is contained in $Sp_\mathcal{B}(s_i)$. To prove the converse inclusion, we start by proving that the spectral radius $R_i\geq 1$ of $s_i$ in $\mathcal{B}$, is actually equal to 1. 

Assume by contradiction $R_i>1$, and let $z_i\in Sp_\mathcal{B}(s_i)$ be such that $|z_i|=R_i$. Then, for $k\geq 1$, the number $z_i^k+z_i^{-k}$ belongs to $Sp_\mathcal{B}(s_i+s_i^{-k})=\linebreak Sp_\mathcal{B}(s_i+(s_i^*)^k)$. Hence:
$$R_i^k-\frac{1}{R_i^k}\leq|z_i^k+z_i^{-k}|\leq \|s_i^k+(s_i^*)^k\|_\mathcal{B},$$
so the norm $\|s_i^k+(s_i^*)^k\|_\mathcal{B}$ grows exponentially in $k$. 
Now by assumption the series $\frac{2}{\pi}(1-\sum_{k=1}^\infty \frac{s_i^k+(s_i^*)^k}{4k^2-1})$ converges in $\mathcal{B}(\Z^n)$, so the norm of the general term has to go to 0. This shows that the norm $\|s_i^k+(s_i^*)^k\|_\mathcal{B}$ grows more slowly than quadratically, which is the desired contradiction.

The same argument applied to $s_i^*$, shows that the spectral radius of $s_i^*=s_i^{-1}$ is also 1, so $Sp_\mathcal{B}(s_i)\subset S^1$.

\item Let $\widehat{\mathcal{B}}$ denote the space of characters of $\mathcal{B}(\Z^n)$, i.e. the set of non-zero, continuous, algebra homomorphisms $\chi:\mathcal{B}(\Z^n)\rightarrow\C$. By density of $\mathcal{B}(\Z^n)$ in $C(\T^n)$, we see that $\T^n\subset\widehat{\mathcal{B}}$, let us see that the converse inclusion also holds. Recall that the {\it joint spectrum} $Sp_\mathcal{B}(s_1,...,s_n)$ of $s_1,...,s_n$ in $\mathcal{B}(\Z^n)$, is:
$$Sp_\mathcal{B}(s_1,...,s_n)=\{(\chi(s_1),...\chi(s_n)):\chi\in\widehat{\mathcal{B}}\}\subset\C^n$$
(see Theorem 3.1.14 in \cite{Horm}). As $\mathcal{B}(\Z^n)$ is topologically generated by $s_1,...,s_n$, the map $\widehat{\mathcal{B}}\rightarrow Sp_\mathcal{B}(s_1,...,s_n):\chi\mapsto(\chi(s_1),...\chi(s_n))$ is a homeomorphism (see Theorem 3.1.15 in \cite{Horm}). Finally, since $\chi(s_i)\in Sp_\mathcal{B}(s_i)$ for every $i$ and every $\chi\in\widehat{\mathcal{B}}$, we have $Sp_\mathcal{B}(s_1,...,s_n)\subset\prod_{i=1}^n Sp_\mathcal{B}(s_i)=\T^n$ by the first part of the proof.

So far we have shown that the inclusion $\mathcal{B}(\Z^n)\hookrightarrow C^*(\Z^n)$ induces an identification of their character spaces. But, in a commutative unital Banach algebra, an element is invertible if and only if no character vanishes on that element (see Theorem 3.1.6 in \cite{Horm}). From this the lemma follows immediately. \epr
\end{itemize}

\begin{Rem} \begin{enumerate}
\item Wiener's celebrated $1/f$-theorem states that $\ell^1(\Z^n)$ is spectral in $C^*(\Z^n)$; since $\ell^1(\Z^n)$ satisfies the assumptions of lemma \ref{spectral}, we see that the latter appears as a generalization of Wiener's $1/f$-theorem\footnote{The reader may have noticed the similarity between the second part of the proof of lemma \ref{spectral}, and Gelfand's proof of Wiener's theorem using commutative Banach algebra techniques (see e.g. Theorem 18.21 in \cite{Rudin}). }
\item Assume that $\mathcal{B}(\Z^n)$ satisfies the assumptions of Proposition \ref{sameKth}; for $n=1,2,3$, the explicit formulae given in this paper directly imply that the map $K_i(\mathcal{B}(\Z^n))\rightarrow K_i(C^*(\Z^n))$ is {\bf onto}. It is trivial for $n=1$. For $n=2$, by formula (\ref{magicproj}) and the assumption on $\mathcal{B}(\Z^2)$, the projection $P_I(s_1,s_2)$ is in $M_2(\mathcal{B}(\Z^2))$. For $n=3$, by Corollary \ref{3torus} we have $$K_0(C^*(\Z^3))=\Z\cdot[1]\oplus\Z\cdot[P_I(s_1,s_2)]\oplus\Z\cdot[P_I(s_2,s_3)]\oplus\Z\cdot[P_I(s_3,s_1)]$$ and $$K_1(C^*(\Z^3))=\Z\cdot[s_1]\oplus\Z\cdot[s_2]\oplus\Z\cdot[s_3]\oplus\Z\cdot[U_I(s_1,s_2,s_3)].$$
By the assumption on $\mathcal{B}(\Z^3)$ together with formulae (\ref{magicproj}) and (\ref{magicunit}), the elements $P_I(s_1,s_2),P_I(s_2,s_3),P_I(s_3,s_1)$ and $U_I(s_1,s_2,s_3)$ are in $M_2(\mathcal{B}(\Z^3))$.




\end{enumerate}\end{Rem}

\begin{Ex} The converse of lemma \ref{spectral} does {\bf not} hold: already for $n=1$, a Banach algebra containing $\C[\Z^n]$, continuously contained in $C^*(\Z^n)$ and spectral in it, does not necessarily contains the root of $2-s_i-s_i^*$. Indeed consider $C^1(S^1)$, the Banach algebra of $C^1$-functions on $S^1$, with the norm $\max\{\|f\|_\infty,\|f'\|_\infty\}$. It is a spectral sub-algebra of $C(S^1)$, as the spectrum of a function is its set of values, both in $C(S^1)$ and $C^1(S^1)$. Since $|\sin{\pi\theta}|$ is not differentiable at $\theta=0$, in $C^1(S^1)$ the element $2-s-s^*$ has no square root \footnote{In Theorem 2.2.10 of \cite{BraRob}, Bratteli and Robinson define the unique positive square root of a positive element $x$ in a unital $C^*$-algebra, as the Riemann integral $x^{1/2}=\frac{1}{\pi}\int_0^\infty x(\lambda+x)^{-1}\frac{d\lambda}{\lambda^{1/2}}$. They don't address the existence of the integral. We observe that $\int_\varepsilon^\infty (\lambda+x)^{-1}\frac{d\lambda}{\lambda^{1/2}}$ exists in any Banach algebra, because of the estimate $\|(\lambda+x)^{-1}\|\leq \frac{1}{\lambda-\|x\|}$ for $\lambda>\|x\|$. So the issue is integrability of $\lambda\ra \frac{(\lambda+x)^{-1}}{\lambda^{1/2}}$ around 0, which holds in a $C^*$-algebra because of the coincidence between norm and spectral radius for positive elements. It is instructive to check why this integrability fails for $x(\theta)=\sin^2(\pi\theta)$ in $C^1[0,1]$.}.
\end{Ex}

\begin{Ex}\label{nonspectral} (J.-B. Bost, see Example 1.1.3 in \cite{JBB}). Fix $\alpha>1$. Let $\ell^1(\Z,w_\alpha)$ denote the weighted $\ell^1$-algebra (a.k.a. Beurling algebra) of $\Z$ corresponding to the weight $w_\alpha(n)=\alpha^{|n|},\;n\in\Z$. Then the inclusion of $\ell^1(\Z,w_\alpha)$ into $C^*(\Z)$ induces isomorphisms in K-theory, by Theorem 1.1.1 in \cite{JBB}. However $\ell^1(\Z,w_\alpha)$ is not spectral in $C^*(\Z)$: indeed, under Fourier transform, $\ell^1(\Z,w_\alpha)$ corresponds to functions on $S^1$ which extend holomorphically to the open corona of radii $\frac{1}{\alpha},\alpha$, and continuously to its closure. It is clear that this algebra is not spectral in $C(S^1)$: for instance $z\mapsto z-\alpha$ is invertible in $C(S^1)$ but not in the subalgebra.
\end{Ex}

\noindent
{\bf Concluding remarks:} The study of embeddings of Banach algebras that induce isomorphisms in K-theory, attracted attention for a while (see Theorems 2.2 and 3.1 in \cite{Swan}; subsection 6.15 in \cite{Kar}; appendix VI.3 in \cite{ConThom}); all those results required the subalgebra to be spectral in one way or the other. The first deep results obtained without any spectral assumption, were obtained by J.-B. Bost \cite{JBB} (see Example \ref{nonspectral} above), in his paper on the Oka principle.

Interest in the subject was revived by V. Lafforgue's seminal work \cite{Laff} on the Baum-Connes conjecture for discrete groups $\Gamma$: he introduced the Banach KK-theory that allowed him, for every unconditional completion $\mathcal{A}(\Gamma)$ of the complex group ring $\C[\Gamma]$, to construct an assembly map $\mu_{\mathcal{A}}:K_i^{top}(\Gamma)\rightarrow K_i(\mathcal{A}(\Gamma))\;(i=0,1)$, where $K_i^{top}(\Gamma)$ is the left-hand side of the Baum-Connes assembly map (see \cite{BCH, GoJuVa}). For $\mathcal{A}(\Gamma)=\ell^1(\Gamma)$, it is a conjecture of J.-B. Bost that $\mu_{\ell^1}$ is an isomorphism for every discrete group $\Gamma$. Note that it is a result by J. Ludwig \cite{Lud} that $\ell^1(\Gamma)$ is spectral in $C^*_r(\Gamma)$ if $\Gamma$ is finitely generated with polynomial growth; and another result by J. Jenkins \cite{Jen} that $\ell^1(\Gamma)$ is NOT spectral in $C^*_r(\Gamma)$ as soon as $\Gamma$ contains the free subsemigroup on 2 generators.

Inspired by all this, Lafforgue introduces a class $\mathcal{C}'$ of groups with geometric interest, for which he proves that $\mu_{\mathcal{A}}$ is an isomorphism for any unconditional completion (Theorem 0.0.2 in \cite{Laff}). To apply this result to the ``official'' Baum-Connes conjecture (see \cite{BCH}), stating that the Baum-Connes assembly map 
$\mu: K_i^{top}(\Gamma)\rightarrow K_i(C^*_r(\Gamma))$ is an isomorphism, he needs to find an unconditional completion that not only embeds into the reduced $C^*$-algebra $C^*_r(\Gamma)$, but has the same K-theory - actually, for groups in the class $\mathcal{C}'$, injectivity of $\mu$ had been known by other means, so only surjectivity of the map $K_i(\mathcal{A}(\Gamma))\rightarrow K_i(C^*_r(\Gamma))$ is relevant. Let us mention lemma 1.7.2 in \cite{Laff}: If $\mathcal{A}(\Gamma)$ is an unconditional completion embedding into $C^*_r(\Gamma)$, a sufficient condition for the map $K_i(\mathcal{A}(\Gamma))\rightarrow K_i(C^*_r(\Gamma))$ to be onto, is that for all $n\geq 1$, every element in $M_n(\C[\Gamma])$ has the same spectral radius in $M_n(\mathcal{A}(\Gamma))$ and in $M_n(C^*_r(\Gamma))$.

Relative spectral conditions like the previous one received some attention in the literature; e.g. the condition that every element in $\C\Gamma$ has the same spectrum in $\ell^1(\Gamma)$ and $C^*_r(\Gamma)$, was studied under the name {\it relatively spectral} by B. Nica \cite{Nica}, and under the name {\it quasi-Hermitian} by H. Samei and M. Wiersma \cite{SaWi}. It is known that such a group is necessarily amenable (Proposition 52 in \cite{Nica}, Theorem 1.4 in \cite{SaWi}); and that finitely generated groups with subexponential growth are completely relatively spectral (i.e. for every $n\geq 1$, any element in $M_n(\C[\Gamma]$ has the same spectrum in $M_n(\ell^1(\Gamma))$ and in $M_n(C^*_r(\Gamma))$); as a consequence the inclusion $\ell^1(\Gamma)\subset C^*_r(\Gamma)$ induces isomorphisms in K-theory (see Example 49 in \cite{Nica}).

\noindent
Authors addresses:\\

\noindent
Gymnase d'Yverdon\\
Route du Gymnase 6\\
CH-1400 Cheseaux-Nor\'eaz\\
Switzerland\\
olivier.isely@eduvaud.ch

\bigskip
\noindent
Institut de math\'ematiques\\
Universit\'e de Neuch\^atel\\
11 Rue Emile Argand - Unimail\\
CH-2000 Neuch\^atel - Switzerland\\
alain.valette@unine.ch


\begin{thebibliography}{CCJJV}

\bibitem[BHS64]{BHS} H. {\sc Bass}, A. {\sc Heller} and R.G. {\sc Swan},
\newblock {\em The Whitehead group of a polynomial extension},
\newblock Publications math\'ematiques de l'I.H.E.S., 22 (1964), 61-79.

\bibitem[BCH94]{BCH}
P. {\sc Baum}, A. {\sc Connes}, and N. {\sc Higson},
\newblock{\em Classifying space for proper actions and K-theory of group $C^*$-algebras},
\newblock Contemp. Math., 167, pp. 240--291, Amer. Math. Soc., Providence, RI, 1994

\bibitem[Bl98]{Black} B. {\sc Blackadar},
\newblock {\em K-Theory for Operator Algebras (2nd ed.)},
\newblock MSRI. Publications 5, Cambridge Univ. Press, 1998.

\bibitem[Bos90]{JBB} J.-B. {\sc Bost}.
\newblock {\em Principe d'Oka, K-th\'eorie et syst\`emes dynamiques non commutatifs},
\newblock Invent. Math. 101 (1990), 261-333.

\bibitem[BR79]{BraRob} O. {\sc Bratteli} and D.W. {\sc Robinson},
\newblock {\em Operator algebras and quantum statistical mechanics I}, 
\newblock Texts and monographs in physics, Springer-Verlag, 1979.



\bibitem[Co81]{ConThom} A. {\sc Connes}, 
\newblock {\em An analogue of the Thom isomorphism for crossed products of a $C^*$-algebra by an action of $\R$},
\newblock Adv. in Math. 39 (1981), 31-55.









\bibitem [GJV19]{GoJuVa} M.P. {\sc Gomez Aparicio}, P. {\sc Julg} and A. {\sc Valette},
\newblock {\em The Baum-Connes conjecture: an extended survey},
\newblock In. A. Chamseddine et al. (eds.), {\it Advances in Noncommutative Geometry}, Springer Nature AG Switzerland 2019, 127-244.



\bibitem[HR00]{HiRo} N. {\sc Higson} and J. {\sc Roe},
\newblock {\em Analytic K-Homology},
\newblock Oxford Math. Monographs, Oxford UP, 2000.

\bibitem[H\"o79]{Horm} L. {\sc H\"ormander},
\newblock {\em An introduction to complex analysis in several variables (2nd ed.)},
\newblock North-Holland Math. Library Vol.7, North-Holland 1979.

\bibitem[Is11]{Isely} O. {\sc Isely},
\newblock {\em K-theory and K-homology for semi-direct products of $\Z^2$ by $\Z$},
\newblock PhD thesis, Universit\'e de Neuch\^atel, 2011.

\bibitem[Je70]{Jen} J. {\sc Jenkins}
\newblock {\em Symmetry and nonsymmetry in the group algebras of discrete groups},
\newblock Pacific J. Math. 32 (1970), 131-145.


\bibitem[Ka78]{Kar} M. {\sc Karoubi},
\newblock {\em K-theory: an Introduction},
\newblock Grundlehren der Math. Wiss., Springer-Verlag, 1978.

\bibitem[La02]{Laff} V. {\sc Lafforgue},
\newblock {\em K-th\'eorie bivariante pour les alg\`ebres de Banach et conjecture de Baum-Connes},
\newblock Invent. Math. 149 (2002), 1-95.

\bibitem[Lo71]{Loday} J.-L. {\sc Loday},
\newblock {\em Applications alg\'ebriques du tore dans la sph\`ere},
\newblock C.R. Acad. Sci., Paris, S\'er. A 272 (1971), 578-581.

\bibitem[Lo88]{Lor} T.A. {\sc Loring},
\newblock {\em K-theory and asymptotically commuting matrices},
\newblock Can. J. Math. XI (1988), 197-216.

\bibitem[Lu79]{Lud} J. {\sc Ludwig},
\newblock {\em A class of symmetric and a class of Wiener group algebras},
\newblock J. Funct. Anal. 31 (1979), 187-194.

 

\bibitem[MSE]{MathStack} Mathematics Stack Exchange,
\newblock {\em Degree 1 map from torus to sphere},
\newblock \verb"https://math.stackexchange.com/questions/1782261/degree-1-"\linebreak \verb"map-from-torus-to-sphere", asked May 12, 2016; retrieved March 7, 2025.

\bibitem[Ni08]{Nica} B. {\sc Nica}
\newblock {\em Relatively spectral morphisms and applications to K-theory},
\newblock Journal of Functional Analysis 255 (2008),s 3303-3328.








\bibitem[PV80]{PV80} M. {\sc Pimsner} and D. {\sc Voiculescu},
\newblock {\em Exact sequences for K-groups and Ext-groups of certain cross-product C*-algebras},
\newblock J. Operator Theory, 4 (1980), 93-118.


\bibitem[Ri81]{Rie} M.A. {\sc Rieffel},
\newblock {\em $C^*$-algebras associated with irrational rotations},
\newblock Pacific J. Math. )1 (1981), 415-429.

\bibitem[Ru66]{Rudin} W. {\sc Rudin},
\newblock {\em Real and complex analysis},
\newblock McGraw-Hill, 1966.

\bibitem[SW20]{SaWi} E. {\sc Samei} and M. {\sc Wiersma},
\newblock {\em Quasi-Hermitian locally compact groups are amenable},
\newblock Advances in Math. 359 (2020), 1-25.

\bibitem[SLL09]{SLL} M.R. {\sc Spiegel}, S. {\sc Lipschutz} and J. {\sc Liu},
\newblock {\em Mathematical handbook of formulas and tables (3rd ed.)}
\newblock Schaum's outlines, McGraw Hill, 2009.




\bibitem[Sw77]{Swan} R. G. {\sc Swan},
\newblock {\em Topological Examples of Projective Modules},
\newblock Trans. Amer. Math. Soc. , 230 (1977), 201-234.








\end{thebibliography}
\end{document}